\documentclass[10pt,a4paper]{article}

\usepackage{graphicx} 
\usepackage{amsmath} 
\usepackage{amssymb}  
\usepackage{arydshln}
\usepackage[affil-it]{authblk}
\usepackage{breakcites}

\newtheorem{lemma}{Lemma}

\begin{document}

\title{On computational issues for stability analysis of LPV systems using parameter dependent Lyapunov functions and LMIs\protect\thanks{This work has been supported by Brazilian funding agencies CAPES and CNPq.}}

\author{Leonardo A. Mozelli\thanks{Email: \texttt{mozelli@cpdee.ufmg.br}; corresponding author}}
\affil{Department of Electronics Engineering,\\ Universidade Federal de Minas Gerais, MG, Brazil}

\author{Ricardo S. L. Adriano}
\affil{Department of Electrical Engineering,\\ Universidade Federal de Minas Gerais, MG, Brazil}

\maketitle

\abstract{This paper deals with the robust stability analysis of linear systems, subject to time-varying parameters. The Parameter Dependent Lyapunov Function are considered, assuming that the temporal derivative of the parameters are bounded. Some computational issues are discussed, which are present in Linear Matrix Inequality (LMI) based approaches and are exacerbated as the quantity of time-varying parameters increases. A possible solution to deal with issues is proposed by modifying the inclusion of the information regarding the time-derivative bounds. Complexity in the number of LMIs constraints can be reduced from very complex to linear. Numerical examples are provide to illustrate the advantages of the proposed methodology.}


\section{INTRODUCTION}\label{sec:intro}

Robust stability analysis plays a central role in systems theory. The theoretical results are useful in several areas: from uncertain systems, whose parameters are invariant but are not precisely known; to linear systems, whose parameters are varying within some limits, as in the case of the Linear Parameter Varying  (LPV) systems; even for nonlinear systems, whose nonlinearities may be embedded by means of appropriate scheduling functions \cite{Shamma:CH:12}. 
To investigate stability, several frameworks were conceived based on obtaining Lyapunov Functions (LFs). A consolidate result is the quadratic stability, which consists on finding a common quadratic LF that guarantees stability for the entire domain of uncertainty, being an important contribution from the 80's \cite{Barmish:JOTA:85}. In the following decade, the procedure of obtaining common LFs to certificate stability was systematized with the advent of optimization tools, giving raise to the Linear Matrix Inequalities (LMIs). 

Over the past decades the LMI framework for robust stability has been intensely researched towards many directions. One direction pointed to the usage of information about time-varying parameters \cite{Gahinet:ITAC:96}. If bounds for the rates of variation of the parameters are available, for instance, less conservative results can be achieved by employing the so called Parameter Dependent Lyapunov Functions (PDLFs). 

PDLFs consist on typical LFs combined using the uncertain parameters. Parametrization can be affine \cite{Gahinet:ITAC:96,Chesi:CDC:04,Geromel:SCL:06,Mozelli:A:09,Han:ITCSRP:14} or polynomial \cite{Chesi:A:07,Montagner:IJC:09}. The shape of LFs used can be quadratic \cite{Gahinet:ITAC:96,Geromel:SCL:06,Mozelli:A:09} or polynomial in the states \cite{Chesi:CDC:04,Chesi:A:07,Montagner:IJC:09,Han:ITCSRP:14}. Recently, some results appear using high-order time-derivatives of the parameters to combine the Lyapunov functions leading to improvements, see \cite{Mozelli:SCL:11}  and references therein.  In \cite{Trofino:IJRNC:2014} general rational dependence on the parameters is considered and in \cite{Pfifer:IJRNC:2015} a griding approach is used to inclued arbitrary dependence on the parameters.

Nevertheless, an aspect rather oversighted is the computational impact of the PDLF for systems with many vertices, or large scale systems. These issues are the main topic of this paper. A simple and scalable example shows that even simple LMI conditions from the Literature can be computationally hard to solve, as they suffer from dimensionality issues. To cope with this effect, a compromise solution is proposed, exploring the geometric structure of the problem. The inclusion of the information from time-derivative of the varying parameters is modified, seeking a balance between conservativeness and computational performance. By this procedure, the numerical complexity is reduced, from a factorial growth to a linear one.  

The reminder of this paper is organized as follows: Section~\ref{sec.pre} gives some theoretical background; some computational issues are discussed in  Section~\ref{sec.issues}; in Section~\ref{sec.sol} a possible solution is proposed and numerical results are presented; finally, section \ref{sec.conclusion} lays down some conclusions and future avenues of research.

\section{BACKGROUND}\label{sec.pre}

Consider uncertain linear systems without input signal:

\begin{equation}
	\dot{x}(t) = A(\theta)x(t),~ \theta\in\mathcal{S}, \label{eq.main_model}
\end{equation}

\noindent with $t\in\mathbb{R}$ being continuous time, $x(t)\in\mathbb{R}^{n}$ the state vector, and the $\theta\in\mathbb{R}^r$ is the vector of uncertainties, belonging to the convex combination:

\begin{equation}
\mathcal{S} = \left\{ s\in\mathbb{R}^n: \sum^r_{i=1}s_i=1,~ s_i\geq 0,~ \forall i=1,2,\ldots,r\right\}. \label{eq.simplex}
\end{equation}

In this paper, uncertainty can be regarded as time-varying, $\theta(t)$, such that this system belongs to the general class of LPV systems \cite{Shamma:CH:12}. Besides, uncertainty is considered in a polytopic fashion. Then function $A: \mathbb{R}^r \rightarrow \mathbb{R}^{n\times n}$ can be modeled by following combination:

\begin{equation}
	A(\theta(t)) = \sum^r_{i=1}\theta_i(t)A_i,
\end{equation}

\noindent for matrices $A_i\in\mathbb{R}^{n\times n}$, $i=1,2,\ldots,r$.

The standard approach to investigate stability of these uncertain systems is to consider a quadratic Lyapunov Function (LF):

\begin{equation}
	V(x(t)) = x^T(t)Px(t) \label{eq.func_std}
\end{equation}

\noindent where $P$ is definite positive, and evaluate if the time-derivative of this LF is always negative, except at the origin:

\begin{equation}
	\dot{V}(x(t)) = x^T(t)\left[A(\theta(t))P+PA(\theta(t))x(t)\right] <0. \label{eq.der_std}
\end{equation}

Since parameters $\theta(t)$ are continuous, constraint \eqref{eq.der_std} has infinite dimension. To circumvent this problem, convexity of the representation can be explored, leading to the following well-known stability condition:

\begin{lemma}\label{lema1}
System \eqref{eq.main_model}  is asymptotically stable if there exists $P = P^T > 0$ satisfying:
\begin{equation}
	A^T_iP+PA_i<0,~\forall i =1,2,\ldots,r. \label{eq.lema1}
\end{equation}
\end{lemma}

\section{SOME COMPUTATIONAL ISSUES}\label{sec.issues}

\subsection{Common Quadratic Lyapunov Function}

The solution presented in Lemma~\ref{lema1} is one of the simplest in Literature. As such, it is very conservative, being unable to guarantee stability for many systems. This happens because no information about the time-varying characteristic of $\theta(t)$ is considered, so a stability certificate is sought even for abrupt (even instantaneous) changes in the uncertain parameters. In this sense, this sort of analysis is strongly related to the field of switched systems, were a common quadratic Lyapunov function is also sought for linear time-invariant sub-systems. In \cite{King:ACC:04} are presented the necessary conditions for the existence of the Lyapunov function for a particular case, for instance.

Even if a common quadratic Lyapunov exists, there are computational issues that  may prevent it from being determined. This is more critical as the number of vertices in the polytope increases too much. The following example, which is scalable both in number of states and of uncertainties, illustrates this dimensionality issue. 

\subsection{Numeric Example}

\begin{figure}
\center
\includegraphics[width=\linewidth]{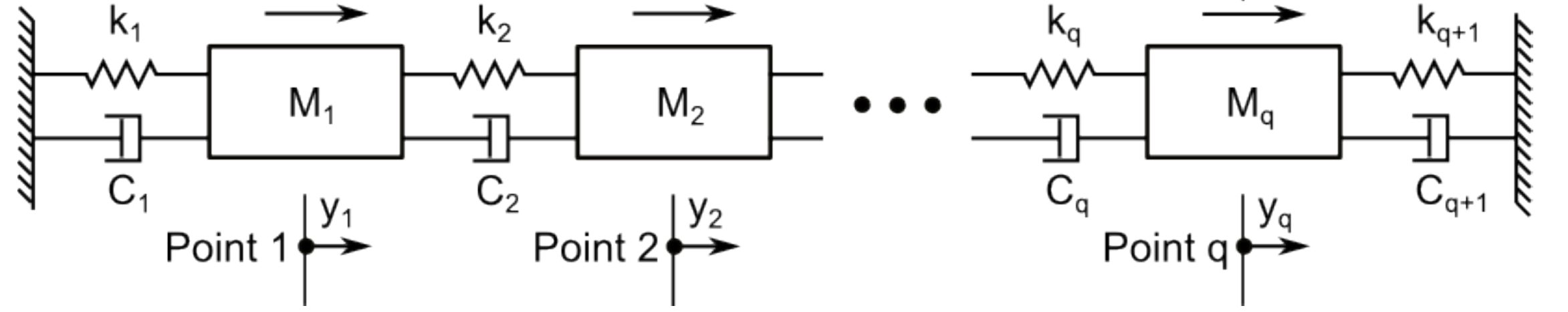}
\caption{Multiple degree of freedom example: mass-damper-spring system.}
\label{fig.car}
\end{figure}

A series of $q$ masses are connected by springs and dampers, as illustrated in Figure~\ref{fig.car}. Each mass weighs $m_j$ and its position is measured over time as $y_j(t)$. The $j$th mass is connected to the next by a spring, with an elastic coefficient  $k_{j+1}$, and by a damper, with a damping coefficient $c_{j+1}$. Special attention must be given to the extreme masses, $m_1$ and $m_q$, which in turn are connected to static structures by dampers $c_1$ and $c_{q+1}$, respectively, and springs $k_1$ and $k_{q+1}$, respectively. This is a basic example of mechanical system with multiple degree of freedom (MDOF). Solving the dynamic free body diagram leads to the equations of motion:

\begin{eqnarray}
	M\ddot{y}(t) + C\dot{y}(t) + Ky(t) = 0, 
\end{eqnarray}  

\noindent where $y(t) = [y_1(t),~y_2(t),~\cdots,~y_q(t)]^T$,

\begin{equation}
 M = 
 \begin{bmatrix} 
 	m_1 & 0 & 0 & \cdots & 0\\
	0 & m_2 & 0 & \cdots & 0\\
	0 & 0 & m_3 & \cdots & 0\\
	\vdots & & &\ddots & \vdots\\ 
	0 & 0 & 0 & \cdots & m_q\\
 \end{bmatrix},
\end{equation}

\begin{equation}
 \Theta(\psi) \!\!=\!\!
 \begin{bmatrix} 
 	\psi_1 +\psi_2 & -\psi_2 & 0 & \cdots & 0\\
	-\psi_2 & \psi_2+\psi_3 & -\psi3 & \cdots & 0\\
	0 & -\psi_3 & \psi_3+\psi_4 & \cdots & 0\\
	\vdots & & &\ddots & \vdots\\ 
	0 & 0 & 0 & \cdots & \psi_{q}+\psi_{q+1}\\
 \end{bmatrix}
\end{equation}

\noindent and $C= \Theta(c)$ and $K= \Theta(k)$.

This system can be rewritten in the state space form as:

\begin{equation}
	\dot{x}(t) =  
	\left[
    	\begin{array}{c;{2pt/2pt}c}
    	0_{n\times n} 	& I_n 
	\\ \hdashline[2pt/2pt]
	-M^{-1}K 		& -M^{-1}C
	\end{array}
    	\right]
  	x(t)
\end{equation}

\noindent with $x(t) = [y(t),~\dot{y}(t)]^T$.

Parameters $k_j$, $m_j$ and $c_j$ can be regarded uncertain within the intervals show in Table~\ref{tab.unc_range}.

\begin{table}[tb]
\caption{Uncertainty range for the parameters in numeric example}
\label{tab.unc_range}
\begin{center}
\begin{tabular}{|l|c|c|c|}
\hline
Parameter & Max & Min & Unit\\
\hline
Elastic Coefficient & 200 & 100 & N/m\\
Damping Coefficient & 8.00 & 4.00 & Ns/m\\
Mass & 5.50 & 5.00 & kg\\
\hline
\end{tabular}
\end{center}
\end{table}

To investigate the impact of dimensionality in the common quadratic analysis, both the number of states and of time-varying parameters are increased. Notice that the number of states increases linearly with the quantity of masses according to $n=2q$. In total, $3q+2$ parameters are necessary to describe a configuration of this system, since the first and last masses are connected by dampers and springs to walls. 

Taking all parameters as time-varying would imply in total of vertices given by $r=2^{3q+2}$ in the system \eqref{eq.main_model}. To provided a more flexible example, only some of the elastic coefficients are considering time-varying. The remaining parameters are considered exact. Therefore, the number os vertices reduces to $r=2^{l}$, where $l\leq q+1$ is the number of parameters that are time-varying.

Bearing this in mind, stability analysis was conducted with Lemma~\ref{lema1}. For each pair $(n,r)$, shown in Table~\ref{tab.feas_quad}, 50 iterations were computed. In each iteration the fixed and time-varying parameters  were chosen in a random fashion within the range shown in Table~\ref{tab.unc_range}. The hardware used was: 2,5GHz \textsc{Intel} Core i5, 4GB 1600MHz DDR, whereas the software was \textsc{Matlab} 2014a, running \textsc{SeDuMi 1.3}  \cite{Sturm:OMS:99}  as solver and  \textsc{Yalmip 2.5}  \cite{Lofberg:CACSD:04}  as parser.

Tables~\ref{tab.feas_quad} and \ref{tab.time_quad} show the feasibility ratio and the average time taken to search for a solution, respectively.

\begin{table}[tb]
\caption{Feasibility ratio: common quadratic stability analysis. Number of states $n$; number of time-varying parameters $l$.}
\label{tab.feas_quad}
\begin{center}
\begin{tabular}{|c|c|c|c|c|}
\hline
$n$ & 4 & 6 & 8 & 10\\
\hline
$r$ &  & & & \\
2 & 1.00 & 1.00 & 1.00 & 1.00\\
4 & 0.98 & 1.00& 1.00 & 0.98\\
8 & 0.94 & 1.00 & 0.98 & 0.88\\
\hline
\end{tabular}
\end{center}
\end{table}

\begin{table}[tb]
\caption{Average computation time (seconds): common quadratic stability analysis. Number of states $n$; number of vertices $r$.}
\label{tab.time_quad}
\begin{center}
\begin{tabular}{|c|c|c|c|c|}
\hline
$n$ & 4 & 6 & 8 & 10\\
\hline
$r$ &  & & & \\
2 &  0.1062 & 0.1105 & 0.1395 & 0.1924\\
4 & 0.1256 & 0.1510 & 0.1906 & 0.2795\\
8 & 0.1784 & 0.2105 & 0.2934 & 0.4985\\
\hline
\end{tabular}
\end{center}
\end{table}

This example reassures the fact that common quadratic LFs are very conservative. It also indicates that increasing  $r$ and $n$, respectively the number of vertices and of states, there is a trend towards greater conservatism. Therefore, dimensionality can be a computational  challenge for this type of LF.

\subsection{Paramenter Dependent Lyapunov Function}

An interesting alternative developed over the years is to consider LFs that depend on the time-varying or uncertain parameters. A possible Parameter Dependent Lyapunov Function (PDLF) is given according to\footnote{In the following time dependency is omitted for sake of clarity}:

\begin{equation}
	V(\theta,x) = 	x^TP(\theta)x = x^T\sum^r_{i=1}\theta_iP_ix
\end{equation}

\noindent which consists on an affine combination of quadratic LFs \cite{Gahinet:ITAC:96}, with the same parametrization used in \eqref{eq.main_model}.

In this case, when calculating the time derivative of the Lyapunov function, information of the time derivative of the parameters appears:

\begin{equation}
	\dot{V}(\theta,x) = x^T\left[A^T(\theta)P(\theta)+P(\theta)A(\theta) +\dot{P}(\theta)\right]x. \label{eq.der_pd}
\end{equation}

This structure has been exploited to improve performance and applied to many nonlinear and time-varying problems \cite{Han:ITCSRP:14,Gaspar:IJC:16,Yang:AST:16}. More recently, high order time-derivatives of the parameters were investigated, producing improvements \cite{Mozelli:SCL:11}. 

From the analytical point of view, the PDLF is more general class and includes the quadratic common as special case. From the computational point of view, it introduces more matrix variables that can relax the analysis. 

\subsection{Time-Derivatives of the Parameters}

Stability analysis  with PDLF when  system \eqref{eq.main_model} posses time-varying parameters is not so straightforward as in the previous sections. Motivation resides in the last term of \eqref{eq.der_pd}: 

\begin{equation}
	\dot{V}_2(\theta,x) = x^T\sum^r_{i=1}\dot{\theta}_iP_ix, \label{eq.dotV2}
\end{equation}

\noindent where 

\begin{equation}
	\sum^r_{i=1}\dot{\theta}_i=0,~ |\dot{\theta_i}|\leq \delta_i,~\forall \theta_i=1,2,\ldots,r, \label{eq.simplex_derivative}
\end{equation}

\noindent and $\delta_i$ are bounds for variations of the uncertain parameters.

Since this condition also has infinite dimension, some approaches to explore its convexity have been proposed over the last years. The more conservative strategy is to consider positive scalar values for the upper bounds, as in \cite{Mozelli:A:09}. However this condition assumes the worst case, with all derivatives being positive, which is not realistic.

A less conservative approach has been proposed in \cite{Chesi:CDC:04,Geromel:SCL:06}. Time-derivatives of the parameters are confined into a manifold with dimension $r-1$, because of the two set of constraints in \eqref{eq.simplex_derivative}: 

\begin{equation}\label{eq.manifold}
    \begin{array}{lcl}
    \Omega &:=& \text{co} \{v^1,v^2,\ldots,v^p\} \\
    \\
    &=& \{v^j \in \mathbb{R}^r |  =\delta_k\leq v^j_k \leq \delta_k, e^Tv^j=0 \},
    \end{array}
\end{equation}

\noindent with $e^T=[1,1,\ldots,1]\in\mathbb{R}^r$, $k$ is the $k$-th coordinate of $v^j_k$. To accumulate this set of vectors  the following matrix is defined:

\begin{equation}\label{eq.Hmatrix}
H :=
\begin{bmatrix}
    v^1, & v^2, & \cdots, &
\begin{pmatrix}
v^j_1\\ v^j_2\\ \vdots\\ v^j_r\\
\end{pmatrix}, & \cdots,
& v^p
\end{bmatrix}\in\mathbb{R}^{r \times p}.
\end{equation}

\vspace*{2mm}

Thus,  combining the columns of matrix $H$ in \eqref{eq.Hmatrix} leads to a finite set of conditions to replace the term \eqref{eq.dotV2}. This became a standard for many researches in the following years. 
 
This polytopic representation has a deep impact over the computational cost: 

\begin{equation}
p = \begin{cases}
\displaystyle \frac{r!}{2\left(\frac{r-1}{2}\right)!},& r \rightarrow \text{is even}\\
\displaystyle \frac{r!}{2\left(\frac{r}{2}\right)!},& r \rightarrow \text{is odd}\\
\end{cases},
\end{equation}

\noindent where $r$ is number of time-varying uncertainties. As Table~\ref{tab.growth} shows, the quantity of columns  $p$ quickly grows to thousands.
 
\begin{table}[tb]
\caption{Vertices growth $p$ with the number of uncertain time-varying parameters $r$.}
\label{tab.growth}
\begin{center}
\begin{tabular}{|c|c|c|c|c|c|c|c|c|c|c|c|}
\hline
$r$ & 2 & 3 & 4 & 5 & 6 & 7 & 8 & 9 & 10 & 11\\
\hline
$p$ & 2 & 6 & 6 & 30 & 20 & 140 & 70 & 630 & 252 & 2772\\
\hline
\end{tabular}
\end{center}
\end{table} 
 
Recently, \cite{Lacerda:IJSS:16} considers a particular case of this general framework. However the conditions are also based on matrix $H$. 
 
The impact of this kind inclusion in the LMIs is better illustrated by the next example. 
 
\subsection{Numeric Example}

The same 600 configurations tested for the standard LF are repeated for the PDLF, using Theorem~1 in \cite{Geromel:SCL:06}, resulting in Tables~\ref{tab.feas_PDLF} and \ref{tab.time_PDLF}, where the feasibility ratio and average computational time are presented, respectively. 

\begin{table}[tb]
\caption{Feasibility ratio: PDLF analysis. Number of states $n$; number of vertices $r$.}
\label{tab.feas_PDLF}
\begin{center}
\begin{tabular}{|c|c|c|c|c|}
\hline
$n$ & 4 & 6 & 8 & 10\\
\hline
$r$ &  & & & \\
2 & 1.00 & 1.00 & 1.00 & 1.00\\
4 & 1.00 & 1.00& 1.00 & 1.00\\
8 & 1.00 & 1.00 & 1.00 & 1.00\\
\hline
\end{tabular}
\end{center}
\end{table}

\begin{table}[tb]
\caption{Average computation time (seconds): PDLF analysis. Number of states $n$; number of vertices $r$.}
\label{tab.time_PDLF}
\begin{center}
\begin{tabular}{|c|c|c|c|c|}
\hline
$n$ & 4 & 6 & 8 & 10\\
\hline
$r$ &  & & & \\
2 &  0.1196 & 0.1453 & 0.1832 & 0.3036\\
4 & 0.3090 & 0.4740 & 0.7799 & 1.4267\\
8 & 10.9121 & 18.2972 & 36.9410 & 106.0802\\
\hline
\end{tabular}
\end{center}
\end{table}

Two aspects are quite noticeable from Tables~\ref{tab.feas_PDLF} and \ref{tab.time_PDLF}. First, the approach using PDLF is indeed less conservative, finding feasible solutions for every parametrization, even in the cases in which the standard LF failed. Secondly, there is a huge impact over the computational cost. In the worst case, the number of vertices in \eqref{eq.simplex_derivative} increased by a factor 35 and the average time taken to solve the stability analysis was increased by a factor of more than 339.

\section{A POSSIBLE SOLUTION}\label{sec.sol}

In the light of the results presented in the previous section, in this section a new approach is proposed to handle the convex inclusion PDLF in terms of LMIs. Since the term in \eqref{eq.dotV2} is responsible for a huge impact in the computational effort, the ideia is to reduce its effect whilst keeping the advantages of the PDLF over the standard quadratic LF.

Instead of considering every vertex in the manifold \eqref{eq.manifold}, the proposed approach consists in obtaining the smaller convex polytope circumscribing the constraints over the time-derivatives of the parameters. In other words, the hypersimplex of dimension $r-1$ circumscribing \eqref{eq.simplex_derivative} is pursued.

For a broarder context, consider that the lower and upper bound of the time-derivative of the parameters can be distinct:

\begin{equation}
	\delta_{i,\min} \leq \dot{\theta}_i \leq \delta_{i,\max}.\label{eq.new_bounds}
\end{equation}

Given constraints in \eqref{eq.simplex_derivative} and in \eqref{eq.new_bounds}, the vertices of the $r-1$ that bounds this manifold are given by:

\begin{equation}
\bar{v}^j = \begin{bmatrix} v^j_{1} & v^j_{2} & \ldots &v^j_{r}\end{bmatrix}^T ,\;\forall j=1,...,r
\end{equation}

\noindent with:

\begin{equation} \label{eq.gen_simplex}
v^j_{i} = \left\lbrace 
\begin{array}{cc}
\triangle_j \cdot \frac{r}{2} + \delta_{j,\min}, & \mbox{if }\; i=j, \\
\\
\delta_{j,\min}, & \mbox{if }\; i\neq j,
\end{array}\right. 
\end{equation}
\noindent where $\triangle_j = \delta_{j,\max} - \delta_{j,\min}$.

Notice that \eqref{eq.gen_simplex} provide a simple way to generate the vertices of the simplex that bounds \eqref{eq.manifold}, not requiring any complex algorithm. These vertices can allocated in columns as in \eqref{eq.Hmatrix}, resulting in:

\begin{equation}\label{eq.barHmatrix}
\bar{H} :=
\begin{bmatrix}
    \bar{v}^1, & \bar{v}^2, & \cdots, &
\begin{pmatrix}
\bar{v}^j_1\\ \bar{v}^j_2\\ \vdots\\ \bar{v}^j_r\\
\end{pmatrix}, & \cdots,
& \bar{v}^p
\end{bmatrix}\in\mathbb{R}^{r \times r}.
\end{equation}

A comparison between \eqref{eq.Hmatrix} and \eqref{eq.barHmatrix} reveals that instead of an exponential increase in the number of columns as in $H$, the proposed approach produces a linear increase in $\bar{H}$.

\subsection{An example for 3D case}

\begin{figure}
\center
\includegraphics[width=0.7\linewidth]{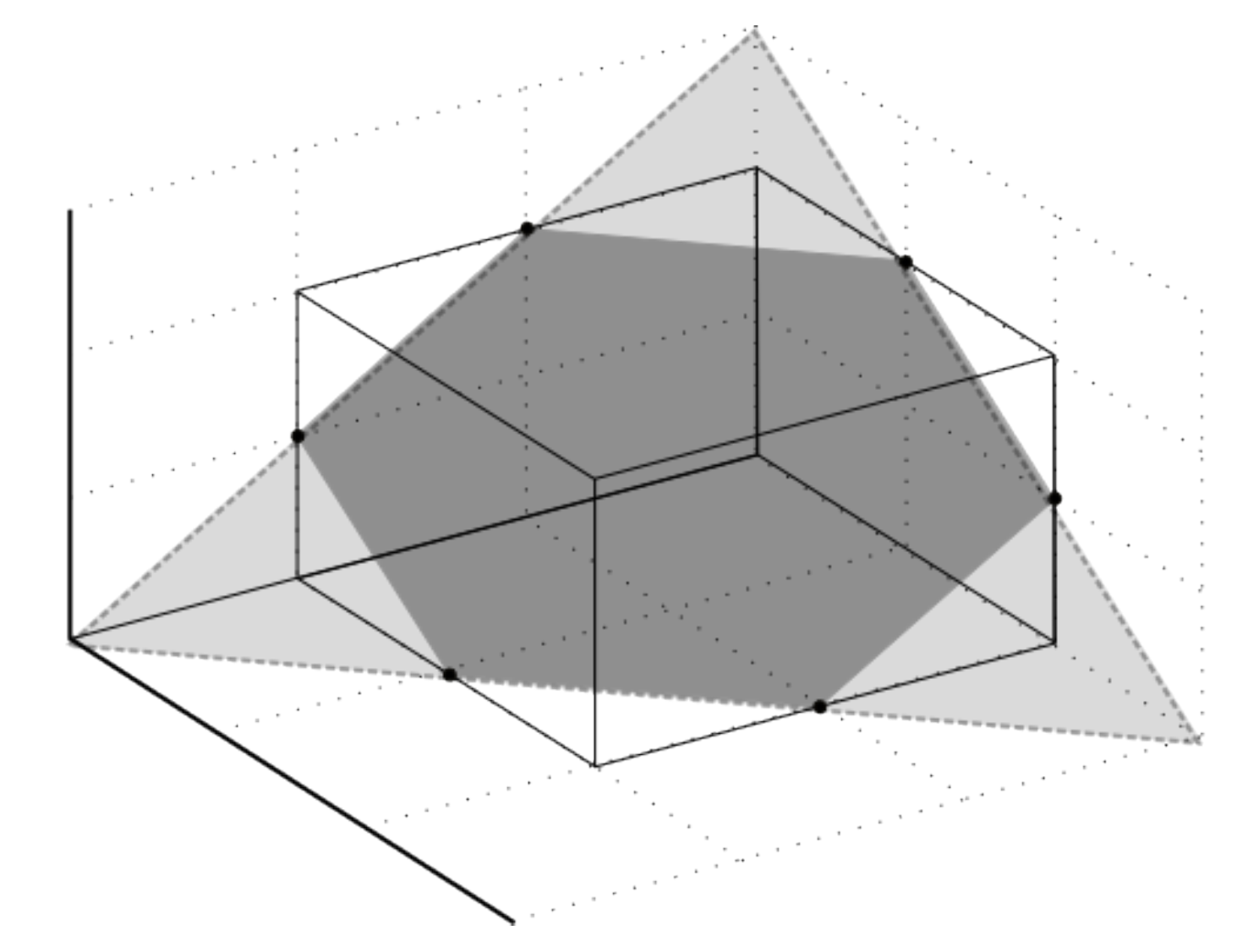}
\caption{Geometric representation of the constraints imposed over the time-derivative of the parameters. In dark grey the manifold associated with $H$; the simplex associated with $\bar{H}$ is in light grey. Notice that overlaps with the dark grey area and contains it.}
\label{fig.3simplex}
\end{figure}

Consider the case of a system \eqref{eq.main_model} with 3 vertices  and $|\dot{\theta}_i| \leq \delta, ~\forall i =1,2,3$.  Constraints \eqref{eq.simplex_derivative} results in the polytope shown in dark grey on Figure~\ref{fig.3simplex}. 

On the other hand, the proposed constraints results in a triangular polytope, which contains the dark grey polytope and the light grey regions. The proposed conditions may be more conservative, since encompasses a larger polytope than the one defined by \eqref{eq.simplex_derivative}. However, the number of vertices grows linearly with $r$, which may be numerically more favorable, as the following example shows.

\subsection{Resuming the Main Numeric Example}

The same 600 configurations tested for the standard LF and PDLF are analyzed a third time, in this turn by the proposed approach. Toward this end, Theorem~1 in \cite{Geromel:SCL:06} has been modified using $\bar{H}$ given by \eqref{eq.barHmatrix} instead of $H$ given by \eqref{eq.Hmatrix}.

The results are listed in Tables~\ref{tab.feas_new} and \ref{tab.time_new}, where the feasibility ratio and average computational time are presented, respectively.

\begin{table}[tb]
\caption{Feasibility ratio: proposed analysis. Number of states $n$; number of vertices $r$.}
\label{tab.feas_new}
\begin{center}
\begin{tabular}{|c|c|c|c|c|}
\hline
$n$ & 4 & 6 & 8 & 10\\
\hline
$r$ &  & & & \\
2 & 1.00 & 1.00 & 1.00 & 1.00\\
4 & 1.00 & 1.00& 1.00 & 1.00\\
8 & 1.00 & 1.00 & 1.00 & 1.00\\
\hline
\end{tabular}
\end{center}
\end{table}

\begin{table}[tb]
\caption{Average computation time (seconds): proposed analysis. Number of states $n$; number of vertices $r$.}
\label{tab.time_new}
\begin{center}
\begin{tabular}{|c|c|c|c|c|}
\hline
$n$ & 4 & 6 & 8 & 10\\
\hline
$r$ &  & & & \\
2 &  0.1173 & 0.1452 & 0.1857 & 0.2983\\
4 & 0.2520 & 0.3382 & 0.5324 & 0.9401\\
8 & 0.7784 & 1.4215 & 2.9696 & 6.0484\\
\hline
\end{tabular}
\end{center}
\end{table}

Both PDLF approaches are able to provide a certificate of stability for all configurations tested in Tables~\ref{tab.feas_PDLF} and \ref{tab.feas_new}. In this example, the performance of both was identical, surpassing the quadratic stability in terms of conservatism, recall Table~\ref{tab.feas_quad}. 

However, major diferences can be spotted when Tables \ref{tab.time_quad}, \ref{tab.time_PDLF} and \ref{tab.time_new} are considered. The proposed approach provides an average time to reach to a solution that is comparable with the one provided by the quadratic LF, although it is alway larger. Yet, the average time is shorter  by
one order of magnitude than the PDLF based on $H$. In the wort case it reduces the average time by a factor of almost 18 times.

Finally, another 50 new configurations using $r=16$ and $n=8$ were tested. For this quantity of uncertainties Theorem~1 with $H$ was unable to finish the computation. However, Theorem~1 with $\bar{H}$ was able to find a feasible solution for every case, with an average computation time of circa $19$ seconds.

\section{CONCLUSIONS}\label{sec.conclusion}

In this paper, a new algorithm has been proposed to numerically describe the constraints related to the temporal derivative of PDLFs. The proposed approach reduced the computational complexity involved with such terms from exponential to linear. The results can be more conservative in some cases, although still are much better than the standard quadratic LF approach. Therefore, it postulates itself as a compromise solution for stability analysis of time-varying linear systems.

Since the advent of PDLF in the context of uncertain LTI systems or LPV systems many results have been produced and applications have benefited.  However, for highly nonlinear, uncertain or complex systems, that require many vertices to be modeled  in the form \eqref{eq.main_model}, some approaches based on PDLF might be prohibitive from the computational point of view. 

The proposed approach tried to balance between the computational efficiency and the reduction of conservatism, presenting a trade-off solution between the quadratic LF and  approaches based on PDLF for time-varying systems. In this way, strategies that already rely on the convex representation give by \eqref{eq.Hmatrix}, as \cite{Han:ITCSRP:14,Gaspar:IJC:16,Yang:AST:16,Mozelli:SCL:11,Mozelli:A:09,Chesi:CDC:04,Geromel:SCL:06,Lacerda:IJSS:16}, besides many others, may benefit from the proposed approach in terms of computational performance for high order systems.

Future lines of research include the extension of the proposed algorithm for  problems other than stability analysis of time-varying linear systems:  control synthesis; computation of performance indexes, such as $H_\infty$ and $H_2$; systems with time-delays or nonlinearities.

\bibliographystyle{apalike}
\bibliography{mozelli.bib}

\end{document}